\ifx\LocalElevenPtArticleFormatLoaded\undefined
  \documentclass[a4paper,11pt]{article}     
  \usepackage{array}                        
  \usepackage{theorem}                      
  \usepackage{amsmath,amscd,amssymb}        
  \usepackage{latexsym}                     
  \usepackage[german,english]{babel}        
  \usepackage[applemac]{inputenc}
 \usepackage{xypic}
\fi

\newtheorem{theorem}{Theorem}[section]
\usepackage{theorem}

\theorembodyfont{\rm}
\newtheorem{remark}{Remark}
\newtheorem{proof}{Proof}

\newcommand{\CC}{{\mathbb{C}}}

\newcommand{\QQ}{{\mathbb{Q}}}

\newcommand{\ZZ}{{\mathbb{Z}}}

\newcommand{\calA}{{\cal A}}

\newcommand{\calM}{{\cal M}}

\newcommand{\op}{\operatorname}
  \newcommand\Igu{{\op{Igu}}}
  \newcommand\Vor{{\op{Vor}}}

\newenvironment{Proof*}[1]{\begin{ProofwCaption}{{#1}}}{\end{ProofwCaption}}
\newenvironment{ProofwCaption}[1]%
  {\addvspace\theorempreskipamount \noindent{\it #1.}\rm}%
  {\qed \par \addvspace\theorempostskipamount}
\newcommand{\qedsymbol}{\mbox{$\Box$}}
\newcommand{\qed}{\quad\qedsymbol}
\begin{document}

\title{Intersection theory of toroidal compactifications of $\calA_4$}
\author{ C. Erdenberger, S. Grushevsky, and K. Hulek}
\date{}
\maketitle

\begin{abstract}
We determine the intersection theory on the Igusa compactification
and the second Voronoi compactification of $\calA_4$.
\end{abstract}

\section{Introduction}

Let $\calA_4$ denote the moduli space of principally polarized
abelian varieties of dimension 4 over $\CC$. There are two natural
toroidal compactifications of $\calA_4$, namely the Igusa
compactification $\calA^\Igu_4$ and the second Voronoi
compactification $\calA^\Vor_4$. The first is given by the perfect
cone or first Voronoi decomposition, which for $g=4$ coincides
with the central cone decomposition, and the latter is given by
the second Voronoi decomposition.

It is well known that
$$
  \op{Pic}(\calA^\Igu_4)\otimes\QQ=\QQ L\oplus\QQ D^\Igu_4
$$
where $D^\Igu_4$ is the closure of the locus rank $1$
degenerations (i.e. semiabelian varieties with the rank of the
abelian part equal to $3$), and $L$ is the $\QQ$-line bundle of
weight $1$ modular forms.

The fan given by the perfect cone decomposition is not basic for
$g=4$, and the second Voronoi decomposition is a basic refinement
of the perfect cone decomposition (this is special to the case
$g=4$). Hence we have a morphism
$$
  \pi : \calA^\Vor_4 \rightarrow \calA^\Igu_4.
$$
One knows (cf. \cite[Proposition 1.6]{HS}) that
$$
  \op{Pic}(\calA^\Vor_4)\otimes\QQ=\QQ L\oplus \QQ
  D^\Vor_4\oplus\QQ E
$$
where $D^\Vor_4$ is again the closure of the rank $1$
degenerations, and $E$ is the divisor contracted to a point under
the map $\pi$. One also knows from \cite[Prop. 1.3]{HS} that
$$
  \pi^{\ast}D^\Igu_4= D^\Vor_4+4E.
$$

Note that the Torelli map gives a morphism
$$
  \overline{t}^\Vor_4: \overline{\cal M}_{4}\rightarrow
  \calA^\Vor_4,
$$
and since we have the morphism $\pi : \calA^\Vor_4 \rightarrow
\calA^ \Igu_4$, we can also define the morphism
$\overline{t}^\Igu_4: \overline{\cal M}_{4}\rightarrow
\calA^\Igu_4$. We denote the images of these maps by ${\cal
J}^\Vor_4$ and ${\cal J}^\Igu_4$ respectively. According to
\cite[Theorem 1.1]{HH}, the classes of these Jacobian loci are
given by
$$
  \left[{\cal J}^\Igu_4 \right] = 8L-D^\Igu_4
$$
and
$$
  \left[{\cal J}^\Vor_4 \right] = 8L-D^\Vor_4 - 4E.
$$
Using these relations and a toric computation we shall determine
the intersection theory of $\calA^\Igu_4$ and $\calA^\Vor_4$. For
genus $4$ this answers a question of N. Shepherd-Barron who posed
in \cite{SB} the problem of determining the intersection theory on
$\calA_g^\Igu$, which he showed to be the canonical model for
$\calA_g$ for $g\ge 12$ (the Picard group of $\calA_g^\Igu$ is
generated by $L$ and $D_g^\Igu$ for all $g\geq 2$).

In the cases $g=2,3$ the Igusa compactification coincides with the
Voronoi compactification and the computation of the intersection
theory can easily be reduced via the Torelli map to calculations
on $\overline{{\cal M}}_g$. In fact for $g=2,3$ the Chow ring of
$\calA_g^\Igu$ is known by the work of van der Geer (\cite{vdG1};
see also \cite{Mu}, \cite{Ts}). We will investigate the
intersection theory of the Igusa compactification $\calA_g^\Igu$
for arbitrary genus in a separate forthcoming paper.

\section{The Igusa compactification}

We shall first treat the Igusa compactification. Let
$$
  a_k:=\langle L^k (D_4^\Igu)^{10-k}\rangle_{\calA_4^\Igu}.
$$
\begin{theorem}
The intersection theory on $\calA_4^\Igu$ is given by\\

\begin{tabular}{|c|c|c|c|c|c|c|c|}
\hline $a_{10}$&$a_9, a_8, a_7$&$a_6$& $a_5, a_4$ &$ a_3$& $a_2$&$a_1$&$a_0$\\
\hline $\vphantom{\dfrac{1}{2}}
\frac{1}{907200}$&$0$&$-\frac{1}{3780}$&$0$&$-\frac{1759}{1680}$&$0$&$\frac{1636249}{1080}$&
$\frac{101449217}{1440}$\\
\hline
\end{tabular}
\end{theorem}
\begin{proof}
We first recall the computation of $a_{10}=L^{10}$. It is well
known (cf. \cite{vdG2},\cite{Ts}) that this number can be computed
using the Hirzebruch-Mumford proportionality theorem (for any
$g$). If $S_k(\Gamma_g)$ denotes the space of weight $k$ cusp
forms, then by \cite[Prop. 2.1]{Ta}
$$
  \op{dim}(S_k(\Gamma_g))\sim 2^{\frac{1}{2}(g-1)(g-2)}\prod^{g}_{j=1}
  \frac{(j-1)!}{(2j)!}B_{2j} \, k^{\frac{1}{2}g(g+1)}
$$
where $B_{2j}$ are the Bernoulli numbers. Hence
$$
  L^{\frac{1}{2}g(g+1)}= \left(\tfrac{1}{2}g(g+1) \right)! \,
  2^{\frac{1}{2}(g-1)(g-2)}\prod^{g}_{j=1}\frac{(j-1)!}{(2j)!} B_{2j}
$$
and a straightforward calculation gives
$$
  a_{10}=L^{10}=\frac{1}{907200}.
$$
We shall now make use of the Torelli map
$$
  \overline{t}_4^\Igu:\overline{{\cal M}}_4\rightarrow \calA_4^\Igu.
$$
We obviously have
$$
  (\overline{t}_4^\Igu)^{\ast}(L)=\lambda, \qquad
  (\overline{t}_4^\Igu)^{\ast} (D_4^\Igu)=\delta_0,
$$
where $\lambda$ is the class of the Hodge bundle on
$\overline{\calM_4}$, and $\delta_0$ is the class of the closure
of the locus of genus $3$ curves with two points identified to
form a node (i.e. the divisor of irreducible stable curves in
$\calM_4$). We let
$$
  b_k:=\langle\lambda^k\delta_0^{g-k}\rangle_{\overline{\calM_4}}.
$$
These intersection numbers on $\overline{\calM_4}$ can be computed
using C. Faber's program implementing his algorithm described in
(\cite{Fa}), and the result is shown in the following table:\\

\noindent
\begin{tabular}{|c|c|c|c|c|c|c|c|c|}
\hline $b_g$ & $b_8, b_7$ & $b_6$ & $b_5$ & $b_4$ & $b_3$ & $b_2$ & $b_1$ & $b_1$ \\
\hline $\vphantom{\dfrac{1}{2}} \frac{1}{113400}$ & $0$&$\frac{1}{3780}$
& $-\frac{2}{945}$& $0$ & $\frac{1759}{1680}$ & $-\frac{1759}{210}$
& $-\frac{1636249}{1080}$ & $-\frac{251987683}{4320}$\\
\hline
\end{tabular}\\

\noindent
Intersecting the relation
$$
  \left[ {\cal J}_4^\Igu\right]= 8L-D_4^\Igu
$$
with $L^{k-1}(D_4^\Igu)^{10-k}$ gives the recurrence relation
$$
  b_{k-1}=8a_k-a_{k-1}.
$$
Since  $a_{10}$ and the numbers $b_{k-1}$ are known, we can
compute the numbers ${a_{k-1}}$ recursively and arrive thus at the
numbers stated.\hfill \qed
\end{proof}
\begin{remark}
In the literature one sometimes finds half the above number for
$L^{\frac{1}{2}g(g+1)}$ --- see \cite{vdG2}. This is the \lq\lq
stack intersection number\rq\rq,  since -$\op{id}$ acts trivially
on the variety $\calA_g$, but defines an involution on the stack.
\end{remark}
\begin{remark}
One can easily see geometrically that $a_9=a_8=a_7=0$. This
follows since $L$ is a pullback of the corresponding line bundle
on the Satake compactification $\calA_g^{\op{Sat}}$, for which the
boundary is codimension 4, and where $L^k$ for $k\geq 7$ can be
represented by a $\QQ$-cycle that does not meet the boundary.
\end{remark}

\section{The Voronoi compactification}

Instead of working with the basis $L, D_4^\Vor$ and $E$ of
$\op{Pic}(\calA_4^\Vor)$, it is easier to use $L, F=
\pi^{\ast}D_4^\Igu$ and $E$. We already saw that
$\pi^{\ast}D_4^\Igu=D_4^\Vor+4E$. Let
$$
  a_{k,l}:=\langle L^k E^l F^{10-k-l}\rangle_{\calA_4^\Vor}.
$$
\begin{theorem}
The intersection theory on $\calA_4^\Vor$ is given by
$$
  a_{k,0}=a_k, \quad a_{k,l}=0 \mbox{ for } 1\leq l\leq 9
$$
$$
  a_{0,10}=E^{10}=-\frac{1680}{1152}=-\frac{35}{24}.
$$
\end{theorem}
\begin{proof}
Since $L$ and $F$ are the pullbacks under $\pi^{\ast}$ of $L$ and
$D_4^\Igu$, we clearly have $a_{k,0}=a_k$. Moreover, $E$ is
contracted to a point in $\calA_4^\Igu$ under $\pi$, and hence
$E^l\cdot\pi^{\ast}(L^k(D_4^\Igu)^{10-k-l})=0$ for all $l$ between
$1$ and $9$.

It thus remains to compute $E^{10}$. This can be done by the
following toroidal computation. The second Voronoi decomposition
for $\calA_4$ is described in \cite{V1}, \cite{V2a}, \cite{V2b},
\cite{ER2}. It is a refinement of the first Voronoi decomposition
obtained by adding another ray $\eta$, generated by the sum of the
primitive generators of the second perfect cone $\Pi_2(4)$. One
then adds all the cones which arise as the span of $\eta$ with
each of the $9$--dimensional faces of $\Pi_2(4)$ and the faces of
these cones, together with their $\mathrm{GL}(4,\ZZ)$ translates.
Thus $\Pi_2(4)$ which was spanned by $12$ rays $\gamma_1, \dots ,
\gamma_{12}$, is divided into sixty-four $10$-dimensional basic
cones. We abuse notations and denote the divisors in the toric
variety defined by this decomposition corresponding to the rays
$\eta$ and $\gamma_1, \dots , \gamma_{12}$ by $E$ and $D_1, \dots
D_{12}$ respectively. Among these divisors we have the following
linear equivalences:
\begin{equation} \label{linrel}
  e_j E + d_{1 \, j} D_1 + \dots + d_{12 \, j} D_{12}
  \sim 0 \quad (j=1, \dots , 10) \;
\end{equation}
where $(e_1, \dots e_{10})$ and $(d_{i \, 1}, \dots , d_{i \,
10})$ are the primitive generators of $\eta$ and $\gamma_i$
respectively (cf. \cite[Theorem (2.3) and section 2.4]{Ts}). We
now multiply equation (\ref{linrel}) successively by $E^9$, then
by $E^8D_i$, $E^7D_iD_j$ with $i \neq j$ up to $ED_{i_1} \cdots
D_{i_8}$. This recursively gives a system of linear equations.
Note that $10$ distinct divisors either intersect transversally in
one point if the $10$ corresponding cones lie in a common
$10$-dimensional cone, or have empty intersection otherwise. This
allows us to solve the above system of linear equations. This can
be done by computer, and gives $E^{10}=-1680$. To get from this
intersection number on the toric variety to the one on the moduli
space, we have to divide by the order of the stabilizer of $\eta$
in $\mathrm{GL}(4,\ZZ)$. The latter is the reflection group $F_4$,
the order of which is equal to 1152 (see Prop. 2.2 of \cite{ER2}
and also the proof of \cite[Prop. 3.6]{HS}.) \hfill \qed
\end{proof}

\bibliographystyle{alpha}

\begin{thebibliography}{HKW}

\bibitem[ER1]{ER1} R.~M.~Erdahl and S.S.~Ryshkov,
The empty sphere, Can. J. Math. {\bf XXXIX} (1987), 794--824.

\bibitem[ER2]{ER2} R.~M.~Erdahl and S.S.~Ryshkov,
The empty sphere, part II, Can. J. Math. {\bf XL} (1988),
1058--1073.

\bibitem[Fa]{Fa} C.~Faber,
Algorithms for computing intersection numbers on moduli spaces of
curves, with an application to the class of the locus of
Jacobians. New trends in algebraic geometry (Warwick, 1996),
93--109, London Math. Soc. Lecture Note Ser., {\bf 264}, Cambridge
Univ. Press, Cambridge, 1999.

\bibitem[vdG1]{vdG1} G.~van der Geer,
The Chow ring of the moduli space of abelian threefolds. J.
Algebraic~Geom.~\textbf{7} (1998), 753--770.

\bibitem[vdG2]{vdG2} G.~van der Geer,
Cycles on the moduli space of abelian varieties. Moduli of curves
and abelian varieties, 65--89, Aspects Math., E33, Vieweg,
Braunschweig, 1999.

\bibitem[HH]{HH} J.~Harris, K.~Hulek,
A Remark on the Schottky Locus in Genus $4$. Proceedings of the
Fano Conference, eds. A.~Colllino, A.~Conte and M.~Marchiso
(2004), 479--483.

\bibitem[HS]{HS} K.~Hulek, G.~K.~Sankaran,
The nef cone of toroidal compactifications of $\calA_4$. Proc.
London Math. Soc. {\bf 88} (2004), 659--704.

\bibitem[Mu]{Mu} D.~Mumford,
Towards an enumerative geometry on the moduli space of curves. In:
Arithmetic and Geometry, vol.~II. Birkh\"auser. Progress in
Mathematics \textbf{36} (1983), 271--328.

\bibitem[S-B]{SB} N.~Shepherd-Barron
Perfect forms and the moduli space of abelian varieties. preprint
math.AG/0502362

\bibitem[Ta]{Ta} Y.-S.~Tai,
On the Kodaira dimension of the moduli space of abelian varieties.
Invent. Math.~\textbf{68} (1982), 425--439.

\bibitem[Ts]{Ts} R.~Tsushima,
A formula for the dimension of spaces of Siegel modular cusp
forms. Am. J. Math.~\textbf{102} (1980), 937--97.

\bibitem[V1]{V1} G.F.~Voronoi,
Nouvelles applications des param\`etres continus \`a la th\'eorie
des formes quadratiques. Premier m\'emoire. Sur quelques
propri\'et\'es des formes quadratiques positives parfaites, J.
Reine Angew. Math. {\bf 133} (1908), 79--178.

\bibitem[V2a]{V2a} G.F.~Voronoi,
Nouvelles applications des param\`etres continus \`a la th\'eorie
des formes quadratiques. Deuxi\`eme m\'emoire. Recherches sur les
parall\'elo\`edres primitifs, J. Reine Angew. Math. {\bf 134},
(1908) 198--287.

\bibitem[V2b]{V2b} G.F.~Voronoi,
Nouvelles applications des param\`etres continus \`a la th\'eorie
des formes quadratiques. Deuxi\`eme m\'emoire. Recherches sur les
parall\'elo\`edres primitifs. Seconde partie. Domaines de formes
quadratiques correspondant aux diff\'erents types de
parall\'elo\`edres primitifs, J. Reine Angew. Math. {\bf 136}
(1909), 67--178.

\end{thebibliography}

\medskip

\noindent
Cord Erdenberger, Klaus Hulek,\\
Institut f\"ur Mathematik (C),\\
Universit\"at Hannover\\
Welfengarten 1, 30060 Hannover, Germany\\
{\tt erdenber@math.uni-hannover.de}, {\tt hulek@math.uni-hannover.de}\\

\medskip

\noindent
Samuel Grushevsky,\\
Mathematics Department,\\
Princeton University,\\
Fine Hall, Washington Road,\\
Princeton, NJ 08544, USA \\
{\tt sam@math.princeton.edu}

\end{document}